\documentclass[12pt]{article}
\usepackage{amssymb, amsmath, amsfonts, tikz, geometry, ytableau}
\usepackage[indent,margin=1cm]{caption}
\usepackage{fullpage}
\usepackage[colorlinks,linkcolor=blue,anchorcolor=blue,citecolor=blue]{hyperref}

\numberwithin{equation}{section}
\numberwithin{figure}{section}

\hypersetup{colorlinks = true}
\newtheorem{thm}{Theorem}[section]
\newtheorem{conj}[thm]{Conjecture}

\newtheorem{lem}[thm]{Lemma}

\def\pf{\noindent{\it Proof.} }
\def\qed{\nopagebreak\hfill{\rule{4pt}{7pt}}
	\medbreak}

\allowdisplaybreaks

\def\pf{\noindent{\it Proof.} }
\def\qed{\nopagebreak\hfill{\rule{4pt}{7pt}}
\medbreak}

\allowdisplaybreaks

\begin{document}
\begin{center}
	{\large \bf  A planar network proof
for Hankel total positivity \\ of type $B$ Narayana polynomials}
\end{center}

\begin{center}
Ethan Y.H. Li$^{1}$, Grace M.X. Li$^{2}$,  Arthur L.B. Yang$^{3}$ and Candice X.T. Zhang$^{4}$\\[6pt]

Center for Combinatorics, LPMC\\
Nankai University, Tianjin 300071, P. R. China\\[8pt]

Email: $^{1}${\tt yinhao\_li@mail.nankai.edu.cn}, $^{2}${\tt limengxing@mail.nankai.edu.cn}, $^3${\tt yang@nankai.edu.cn}, $^{4}${\tt zhang\_xutong@mail.nankai.edu.cn}
\end{center}

\noindent\textbf{Abstract.}
The Hankel matrix of type $B$ Narayana polynomials was proved to be totally positive by Wang and Zhu, and independently by Sokal. Pan and Zeng raised the problem of giving a planar network proof of this result. In this paper, we present such a proof by constructing a planar network allowing negative weights, applying the Lindstr\"{o}m-Gessel-Viennot lemma and establishing an involution on the set of nonintersecting families of directed paths.



\section{Introduction}\label{intro}

For any $n\geq 0$, let
\begin{align*}
W_n(q)=\sum_{k=0}^{n}\binom{n}{k}^2 q^k
\end{align*}
denote the $n$-th Narayana polynomial of type $B$. Wang and Zhu \cite{WZ16}, and Sokal \cite{Sok15} independently proved that
the Hankel matrix \begin{equation}\label{eq-Hankel-def}
H =(W_{i+j}(q))_{i,j \ge 0}
\end{equation}
is $q$-totally positive, namely, any minor of $H$ is a polynomial in $q$ with nonnegative coefficients.
The main objective of this paper is to give a combinatorial proof of the $q$-total positivity of $H$, which solves a problem of Pan and Zeng \cite{PZ16}.

The $q$-total positivity of the Hankel matrix $H$ arose in the study of the $q$-log-convexity of the polynomial sequence $(W_n(q))_{n\geq 0}$. For the convenience of introducing related definitions and results, we make use of the notion of $q$-nonnegativity and the symbol $\ge_q$. A polynomial $f(q)$ with real coefficients is called $q$-nonnegative, written $f(q) \ge_q 0$, if all its coefficients are nonnegative. Accordingly, for two polynomials $f(q)$ and $g(q)$ we write $f(q) \ge _q g(q)$ if $f(q)-g(q) \ge_q 0$. Recall that a sequence
$\alpha = (a_n(q))_{n\geq 0}$ of polynomials in $q$ is said to be $q$-log-convex if for any $n\geq 1$ there holds
$a_{n+1}(q)a_{n-1}(q)\geq_q a^2_n(q)$. Furthermore, if
$a_{m+1}(q)a_{n-1}(q)\geq_q a_m(q)a_n(q)$ holds for any $m\geq n\geq 1$, then we call $\alpha$ a strongly $q$-log-convex sequence. Conversely,
we say that $\alpha$ is a $q$-log-concave sequence if for any $n\geq 1$ we have
$a^2_n(q)\geq_q a_{n+1}(q)a_{n-1}(q)$, and it is a strongly $q$-log-concave sequence
if $a_m(q)a_n(q)\geq_q a_{m+1}(q)a_{n-1}(q)$ holds
for any $m\geq n\geq 1$. The concept of $q$-log-concavity was introduced by Stanley, and the notion of strong $q$-log-concavity was due to Sagan \cite{Sag92DM}. Many polynomial sequences have been proved to be $q$-log-concave, or even strongly $q$-log-concave, see Butler \cite{But90}, Krattenthaler \cite{Kra89}, Leroux \cite{Ler90}, Sagan \cite{Sag92DM, Sag92TAMS}, and Chen, Wang and Yang \cite{CWY11}. However, $q$-log-convex sequences received very little attention until the work of Liu and Wang \cite{LW07}, who first introduced the notion of $q$-log-convexity. Liu and Wang established the $q$-log-convexity of many combinatorial polynomials, such as the Eulerian polynomials. For further progress on $q$-log-convexity, see \cite{CWY10, Zhu14} for instance.

The $q$-log-convexity of $(W_n(q))_{n\geq 0}$ was conjectured by Liu and Wang \cite{LW07}, and was proved later by Chen, Tang, Wang and Yang \cite{CTWY10} by using the theory of symmetric functions. Zhu \cite{Zhu13} further established the strong $q$-log-convexity of $(W_n(q))_{n\geq 0}$ by identifying this polynomial sequence as the first column of the triangular array
$B=(b_{n,k}(q))_{n,k \ge 0}$, which is generated by
\begin{equation}\label{eq-narab-recurrence}
\begin{split}
b_{n,0}(q)&= (1+q)\cdot b_{n-1,0}(q) + 2q \cdot b_{n-1,1}(q);\\
b_{n,k}(q)&= b_{n-1,k-1}(q) + (1+q)\cdot b_{n-1,k}(q) + q \cdot b_{n-1,k+1}(q) \quad (k\ge 1, \, n\ge 1)
\end{split}
\end{equation}
with $b_{0,0}(q) = 1$ and $b_{n,k}(q) = 0$ for $k > n$. The triangular array $B$ belongs to a wide class of matrices, called $q$-recursive matrices in \cite{WZ16}, or Catalan-Stieltjes matrices in \cite{PZ16, LLYZ21}, which we will recall below.
Let $\gamma=(r_k(q))_{k\ge 0}$, $\sigma=(s_k(q))_{k\ge 0}$ and $\tau=(t_k(q))_{k\ge 1}$ be three sequences of polynomials in $q$. The Catalan-Stieltjes
matrix with respect to $\gamma,\sigma,\tau$, denoted by  $C^{\gamma,\sigma,\tau}=(c_{n,k}(q))_{n,k \ge 0}$, is generated by the following recursive relations:
\begin{equation*}
\begin{split}
c_{n,0}(q) &=s_0(q) c_{n-1,0}(q) + t_1(q)c_{n-1,1}(q);\\
c_{n,k}(q) &= r_{k-1}(q)c_{n-1,k-1}(q) + s_k(q)c_{n-1,k}(q) + t_{k+1}(q)c_{n-1,k+1}(q) \quad (k\ge 1, \, n\ge 1),
\end{split}
\end{equation*}
where $c_{0,0}(q)=1$ and $c_{n,k}(q)=0$ unless $n\geq k\geq 0$. Actually, Zhu \cite{Zhu13} gave a general criterion for the strong $q$-log-convexity of $(c_{n,0}(q))_{n\geq 0}$ of
$C^{\gamma,\sigma,\tau}$. Further, Wang and Zhu \cite{WZ16} proved that the Hankel matrix $(c_{i+j,0}(q))_{i,j\geq 0}$ is $q$-totally positive provided that the matrix
$$
L^{\gamma,\sigma,\tau}=\begin{pmatrix}
1 &\quad &\quad &\quad&\quad\\
s_0(q) & r_0(q) & \quad &\quad &\quad\\
t_1(q) & s_1(q) & r_1(q)  &\quad &\quad \\
\quad& t_2(q) & s_2(q) & r_2(q) &\quad\\
\quad&\quad&\ddots & \ddots &\ddots
\end{pmatrix},
$$
called the coefficient matrix of $C^{\gamma,\sigma,\tau}$,
is $q$-totally positive. As a result, Wang and Zhu  
obtained the $q$-total positivity of the Hankel matrix $H=(W_{i+j}(q))_{i,j\geq 0}$, which was also independently proved by Sokal \cite{Sok15} based on the continued fraction expression of the generating function $\sum_{n\geq 0}W_n(q)x^n$.

We would like to note that the $q$-total positivity of the Hankel matrix $(c_{i+j,0}(q))_{i,j \ge 0}$ is also closely related to that of $C^{\gamma,\sigma,\tau}$, for details see \cite{LMW16} and \cite{WZ16}. Chen, Liang and Wang \cite{CLW15} raised the problem of giving a combinatorial interpretation for the $q$-total positivity of $C^{\gamma,\sigma,\tau}$. An ideal tool to combinatorially proving the positivity of a matrix is the famous Lindstr\"om-Gessel-Viennot lemma, see \cite{Lin73, GV85, GV89}. A natural strategy is to construct a planar network with nonnegative weights for the target matrix, and then to apply the Lindstr\"om-Gessel-Viennot lemma to interpret each minor of this matrix as the generating function of nonintersecting families of directed paths, which are obviously nonnegative. In the spirit of this method, Pan and Zeng \cite{PZ16} provided a general planar network construction for the Catalan-Stieltjes matrices and their associated Hankel matrices, which enables them to give combinatorial proofs of the $q$-total positivity for many such matrices, such as those related to the Eulerian polynomials, Schr\"oder polynomials, and Narayana polynomials of type $A$. However,
their approach did not work for Narayana polynomials of type $B$, and they proposed it as an open problem to find a planar network proof of the $q$-total positivity of $H=(W_{i+j}(q))_{i,j\geq 0}$.
It seems impossible to find a planar network with only nonnegative weights for $H$.

In this paper, inspired by our recent work \cite{LLYZ21}, we construct for $H$ a suitable planar network allowing negative weights and solve Pan and Zeng's problem. 
In our construction, the planar network for $H$ can be naturally divided into serial segments which are essentially subnetworks of the planar network for the coefficient matrix $L^B$ of $B$.
By applying the Lindstr\"om-Gessel-Viennot lemma and establishing a sign-reversing involution on the nonintersecting families of each segment, we combinatorially prove the $q$-total positivity of $L^B$ and $H$.

This paper is organized as follows.
In Section \ref{sect-lgv}, we will introduce the Lindstr\"{o}m-Gessel-Viennot lemma. In Section \ref{sect-network}, we will present our planar network construction for the coefficient matrix $L^B$, as well as a combinatorial proof of its $q$-total positivity. In Section \ref{sect-pf} we will make use of
the results in Section \ref{sect-network}
to obtain a planar network for $H$ and a combinatorial proof of its $q$-total positivity.
We conclude this paper in Section \ref{sect-conj} with a conjecture on the immanant positivity for
$H$.

\section{The Lindstr\"om-Gessel-Viennot lemma}\label{sect-lgv}

The Lindstr\"om-Gessel-Viennot lemma was originally proved by Lindstr\"om \cite{Lin73} and further developed by Gessel and Viennot \cite{GV85, GV89}. It has a broad range of applications, see \cite{HG95, Kra96, MW00, Ste90} for instance. In this section, we will give an overview of the Lindstr\"om-Gessel-Viennot lemma, which plays a key role in our combinatorial proof of the $q$-total positivity of the Hankel matrix of type $B$ Narayana polynomials.

To state the Lindstr\"om-Gessel-Viennot lemma, we need some notations. Let $D$ be a directed graph, or digraph for short, with vertex set $V(D)$ and arc set $A(D)$. A digraph $D$ is said to be acyclic if it contains no directed cycles. Throughout this paper we may assume that $D$ is locally finite, namely, for any two vertices $u,v \in V(D)$ the number of directed paths from $u$ to $v$ is finite. We say two directed paths intersect if they have a vertex in common. A sequence $(p_1,\ldots,p_n)$ of directed paths is called a nonintersecting family if $p_i$ and $p_j$ do not intersect for any $i \neq j$. Let $\mathbf{U} = (u_1,\ldots,u_n)$ and $\mathbf{V} = (v_1,\ldots,v_n)$ be two sequences of vertices in $D$, and let $\mathcal{N}_D(\mathbf{U},\mathbf{V})$ denote the set of nonintersecting families $(p_1,\ldots,p_n)$ such that $p_i$ is a directed path from $u_i$ to $v_i$ for each $1 \le i \le n$. If for any permutation $\sigma$ of $\{1,2,\ldots,n\}$, the set $\mathcal{N}_D(\mathbf{U},\sigma(\mathbf{V})) = \mathcal{N}_D((u_1,\ldots,u_n),(v_{\sigma(1)},\ldots,v_{\sigma(n)}))$ is empty unless $\sigma$ is the identity permutation, then $\mathbf{U}$ and $\mathbf{V}$ are said to be compatible. A weight function $\mathrm{wt}$ of $D$ is a map from $A(D)$ to $R$, where $R$ is a commutative ring with identity. The weight of a directed path in $D$ is the product of the weights of all its arcs, and the weight of a nonintersecting family is defined to be the product of the weights of all its components. Given two vertices $u$ and $v$ of $D$, let $GF_D(u,v)$ denote the sum of the weights of all directed paths from $u$ to $v$.
For two sequences $\mathbf{U}$ and $\mathbf{V}$ of vertices in $D$, let $GF(\mathcal{N}_D(\mathbf{U},\mathbf{V}))$ denote the sum of the weights of all elements in $\mathcal{N}_D(\mathbf{U},\mathbf{V})$. For a matrix $M$, we denote by $\det [M]$ the determinant of $M$. The celebrated Lindstr\"om-Gessel-Viennot lemma is stated as follows.

\begin{lem}[{\cite[Corollary 2]{GV89}}]\label{lem-lgv}
Let $D$ be a locally finite and acyclic digraph with a weight function, and let $\mathbf{U} = (u_1,\ldots,u_n)$, $\mathbf{V} = (v_1,\ldots,v_n)$ be two sequences of vertices in $D$. Then
\[
\det \left[\left(GF_D(u_i,v_j)\right)_{1 \le i,j \le n}\right] = GF(\mathcal{N}_D(\mathbf{U},\mathbf{V})).
\]
\end{lem}

In this paper, we mainly apply the above lemma to a special class of digraphs, called planar networks. Recall that a digraph $D$ is said to be planar if it can be embedded in the plane with edges meeting only at endpoints. We call $\mathcal{D} = (D,\mathrm{wt}_D)$ a planar network if $D$ is a locally finite, acyclic, and planar digraph, and $\mathrm{wt}_D$ is a weight function of $D$. Given an $n \times n$ matrix $M$, then $\mathcal{D}$ is called a planar network for $M$ if there exist two sequences $(u_1,\ldots,u_n)$ and $(v_1,\ldots,v_n)$ of vertices in $D$ such that
\begin{equation*}
M = \left(GF_D(u_i,v_j)\right)_{1 \le i,j \le n}.
\end{equation*}
In the remaining part of this paper, we usually specify the vertices and say that $\mathcal{D} = (D,\mathrm{wt}_D,(u_1,\ldots,u_n),(v_1,\ldots,v_n))$ is a planar network for $M$.

\section{The coefficient matrix \texorpdfstring{$L^B$}{}}\label{sect-network}

In this section we will establish the planar network for the coefficient matrix $L^B$ and prove its $q$-total positivity. By \eqref{eq-narab-recurrence}, we have
\begin{equation*}
L^B=(l_{i,j})_{i,j \ge 0} =
\begin{pmatrix}
1 & & & & &\\
1+q & 1 & & & &\\
2q & 1+q & 1 & & & \\
& q & 1+q & 1 & &\\
& & q & 1+q & 1 & \\
& & &\ddots & \ddots &\ddots\\
\end{pmatrix}.
\end{equation*}

\begin{figure}[htb]
  \centering
\begin{tikzpicture}
[place/.style={thick,fill=black!100,circle,inner sep=0pt,minimum size=1mm,draw=black!100},scale=1.5]
\draw [thick] [->] (2.5,2.5) -- (3,2.5);
\draw [thick] [->] (3,2.5) -- (3.5,2.5) -- (4,2.5);
\draw [thick] (4,2.5) -- (4.5,2.5);
\node [place,label=above:{\footnotesize$P_0$}] at (2.5,2.5) {};
\node [place,label=above:{\footnotesize$Q_0$}] at (3.5,2.5) {};
\node [place,label=above:{\footnotesize$P'_0$}] at (4.5,2.5) {};
\draw [thick] [->](2.5,1.5) -- (3,1.5);
\draw [thick] [->] (3,1.5) -- (4,1.5);
\draw [thick] (4,1.5) -- (4.5,1.5);
\node [place,label=below:{\footnotesize$P_1$}] at (2.5,1.5) {};
\node [place,label=below:{\footnotesize$Q_1$}] at (3.5,1.5) {};
\node [place,label=below:{\footnotesize$P'_1$}] at (4.5,1.5) {};
\draw [thick] [->] (2.5,1.5) --(3,2);
\draw [thick] (3,2) --  (3.5,2.5) ;
\draw [thick] [->] (2.5,1.5) -- (3.5,2);
\draw [thick] (3.5,2) -- (4.5,2.5);
\draw [thick] [->] (2.5,0.5) -- (3,0.5);
\draw [thick] [->] (3,0.5) -- (4,0.5);
\draw [thick] (4,0.5) -- (4.5,0.5);
\node [place,label=below:{\footnotesize$P_2$}] at (2.5,0.5) {};
\node [place,label=below:{\footnotesize$Q_2$}] at (3.5,0.5) {};
\node [place,label=below:{\footnotesize$P'_2$}] at (4.5,0.5) {};
\draw [thick] [->] (2.5,0.5) -- (3,1);
\draw [thick] (3,1) -- (3.5,1.5);
\draw [thick] [->] (3.5,0.5) -- (4,1);
\draw [thick] (4,1) -- (4.5,1.5);
\draw [thick] [->] (2.5,-0.5) -- (3,-0.5);
\draw [thick] [->] (3,-0.5) -- (4,-0.5);
\draw [thick]  (4,-0.5) -- (4.5,-0.5);
\node [place,label=below:{\footnotesize$P_3$}] at (2.5,-0.5) {};
\node [place,label=below:{\footnotesize$Q_3$}] at (3.5,-0.5) {};
\node [place,label=below:{\footnotesize$P'_3$}] at (4.5,-0.5) {};
\draw [thick] [->] (2.5,-0.5) -- (3,0);
\draw [thick] (3,0) -- (3.5,0.5);
\draw [thick] [->] (3.5,-0.5) -- (4,0);
\draw [thick] (4,0) -- (4.5,0.5);
\node [blue] at (2.75,2) {$q$};
\node [blue] at (3.375,2.15) {$-1$};
\node [blue] at (2.75,1) {$q$};
\node [blue] at (2.75,0) {$q$};
\draw [thick]plot[smooth, tension=.7] coordinates {(3.5,1.5) (4.1,1.9) (4.5,2.5)};
\draw [thick]plot[smooth, tension=.7] coordinates {(3.5,1.5) (3.9,2.1) (4.5,2.5)};
\draw [thick][->] (4.1,1.9) -- (4.11,1.91);
\draw [thick][->] (3.9,2.1) -- (3.91,2.11);
\node at (3.5,-1) {$\textbf{\vdots}$};
\node at (2.5,-1) {$\textbf{\vdots}$};
\node at (4.5,-1) {$\textbf{\vdots}$};
\end{tikzpicture}
\caption{The planar network $\mathcal{L^B}$}\label{fig-network-lb}
\end{figure}
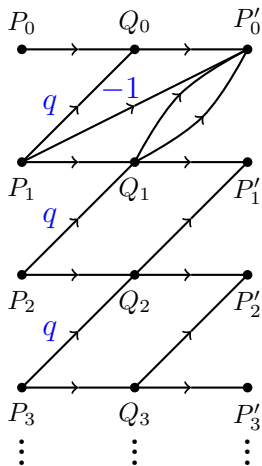

Now we give the construction of the planar network for $L^B$. Let $D^{L^B}$ be the infinite planar digraph with vertex set
\[
V(D^{L^B})=\{P_i \mid i \ge 0\} \cup \{Q_i \mid i \ge 0\} \cup \{P'_i \mid i \ge 0\}
\]
and arc set
\begin{align*}
A(D^{L^B})=&
\{P_i \to Q_i \mid i \ge 0\} \cup \{P_i \to Q_{i-1} \mid i \ge 1\} \\
& \cup \{Q_i \to P'_i \mid i \ge 0\} \cup \{Q_i \to P'_{i-1} \mid i \ge 2\} \\
& \cup \{P_1 \to P'_0, Q_1 \overset{l}{\to} P'_0, Q_1 \overset{r}{\to} P'_0\},
\end{align*}
where the coordinates of vertices are given by $P_i = (0,-i)$, $Q_i = (1,-i)$ and $P'_i = (2,-i)$, and $Q_1 \overset{l}{\to} P'_0$, $Q_1 \overset{r}{\to} P'_0$ are multiple arcs from $Q_1$ to $P'_0$ with one drawn on the left and the other on the right, respectively, as shown in Figure \ref{fig-network-lb}. The weight function $\mathrm{wt}_{D^{L^B}}$ is defined by
\begin{align*}
\mathrm{wt}_{D^{L^B}} (P_1 \to P'_0)=-1, \quad \mathrm{wt}_{D^{L^B}} (P_i \to Q_{i-1})= q \quad \mathrm{for} \quad i \ge 1,
\end{align*}
and $\mathrm{wt}_{D^{L^B}}(a) = 1$ for the other arcs $a$ in $A(D^{L^B})$. Then we have the following result.
\begin{lem}\label{lem-network-lb}
Let $D^{L^B}$ and $\mathrm{wt}_{D^{L^B}}$ be defined as above. Then
\[
\mathcal{L^B} = (D^{L^B},\mathrm{wt}_{D^{L^B}},(P_0,P_1,\ldots),(P'_0,P'_1,\ldots))
\]
is a planar network for $L^B$, or equivalently,
\[
L^B = \left(GF_{D^{L^B}}(P_i,P'_j)\right)_{i,j \ge 0}.
\]
\end{lem}
\pf By the definitions in Section \ref{sect-lgv} and the above construction, it is straightforward to verify that
\[
l_{i,j} = GF_{D^{L^B}}(P_i,P'_j)
\]
for $i,j \ge 0$. Then the proof follows. \qed

Figure \ref{fig-network-lb} provides an illustration of the planar network $\mathcal{L^B}$, where we only label the weights not equal to 1.

The remaining part of this section is devoted to giving a combinatorial proof of the $q$-total positivity of $L^B$ by using the planar network $\mathcal{L^B}$. To this end, we define the following three properties of nonintersecting families $\mathbf{p}=(p_1,\ldots,p_k)$:
\begin{itemize}
\item[($\mathcal{P}_1$)] $p_1 = P_1\to P'_0$;
\item[($\mathcal{P}_2$)] $p_1 = P_1 \to Q_1 \overset{l}{\to} P'_0$;
\item[($\mathcal{P}_3$)] there exists $l$ ($l \ge 2$) such that $p_m = P_m \to Q_{m-1} \to P'_{m-1}$ for $1 \le m \le l-1$ and $p_l = P_l \to Q_l \to P'_{l-1}$.
\end{itemize}
Given a positive integer $k$ and two sequences
\begin{equation}\label{eq-ij}
\begin{split}
I& = (i_1,\ldots,i_k),\, \mbox{ where } 0 \le i_1 < \cdots < i_k,\\
J& = (j_1,\ldots,j_k),\, \mbox{ where } 0 \le j_1 < \cdots < j_k,
\end{split}
\end{equation}
let
\begin{align}\label{eq-pipjp}
\mathbf{P}_I=(P_{i_1},\ldots,P_{i_k}),\quad
\mathbf{P}'_J=(P'_{j_1},\ldots,P'_{j_k}),
\end{align}
and let
\begin{align}\label{eq-sij}
S_{I,J}& = \{\mathbf{p} \in \mathcal{N}_{D^{L^B}}(\mathbf{P}_I,\mathbf{P}'_J) \mid \mathbf{p} \mbox{ satisfies none of ($\mathcal{P}_1$), ($\mathcal{P}_2$), ($\mathcal{P}_3$)}\}.
\end{align}
It is clear that each nonintersecting family in $S_{I,J}$ has a $q$-nonnegative weight. By virtue of this, the following result provides a combinatorial interpretation for the $q$-total positivity of $L^B$.

\begin{thm}\label{thm-qtp-lb}
Let $I,J, \mathbf{P}_I,\mathbf{P}'_J, S_{I,J}$ be as given in \eqref{eq-ij}, \eqref{eq-pipjp} and \eqref{eq-sij} respectively, and let $L^B_{I,J}$ denote the submatrix of $L^B$ whose rows are indexed by $I$ and columns indexed by $J$. Then we have
\begin{align}\label{eq-lbij-det}
\det \left[L^B_{I,J}\right] = GF(S_{I,J}),
\end{align}
where $GF(S_{I,J})$ denotes the sum of weights of all elements in $S_{I,J}$. In particular, $\det \left[L^B_{I,J}\right]$ is $q$-nonnegative.
\end{thm}

\pf By Lemma \ref{lem-network-lb} and Lemma \ref{lem-lgv}, we have
\[
\det \left[L^B_{I,J}\right] = GF(\mathcal{N}(\mathbf{P}_I,\mathbf{P}'_J)),
\]
where we use
$\mathcal{N}(\mathbf{P}_I,\mathbf{P}'_J)$ to stand for $\mathcal{N}_{D^{L^B}}(\mathbf{P}_I,\mathbf{P}'_J)$ for convenience.
In the following we may assume that $\mathcal{N}(\mathbf{P}_I,\mathbf{P}'_J)\neq \emptyset$, otherwise, $S_{I,J}=\emptyset$, $\det \left[L^B_{I,J}\right] = 0$, and \eqref{eq-lbij-det} holds trivially.
For $i=1,2,3$, let
\begin{align*}
\mathcal{N}_i&=\{\mathbf{p} \in \mathcal{N}(\mathbf{P}_I,\mathbf{P}'_J) \mid \mathbf{p} \mbox{ satisfies } (\mathcal{P}_i)\}.
\end{align*}
Clearly, $\mathcal{N}(\mathbf{P}_I,\mathbf{P}'_J)$ is the disjoint union of
$\mathcal{N}_1,\mathcal{N}_2,\mathcal{N}_3$ and $S_{I,J}$.
Now it suffices to give an involution $\phi$ on $\mathcal{N}(\mathbf{P}_I,\mathbf{P}'_J)$ such that
any nonintersecting family $\mathbf{p}\in \mathcal{N}_1\cup \mathcal{N}_2\cup\mathcal{N}_3$ and its image $\phi(\mathbf{p})$
have opposite weights, and the restriction of $\phi$ to $S_{I,J}$ is the identity map. Thus, we only need to define the action of $\phi$ on $\mathcal{N}_1\cup \mathcal{N}_2\cup\mathcal{N}_3$.

Let us first consider the case $k=1$, for which we have $\mathcal{N}_3=\emptyset$. If $i_1=1,j_1=0$, then we have
$$\mathcal{N}_1=\{P_1\to P_0'\},\quad \mathcal{N}_2=\{P_1\to Q_1 \overset{l}{\to} P'_0\}.$$ Define $\phi$ to be the map which sends $P_1\to P_0'$ and $P_1\to Q_1 \overset{l}{\to} P'_0$ to each other. Note that the weight of $P_1\to P_0'$ is $-1$, while the weight of $P_1\to Q_1 \overset{l}{\to} P'_0$ is $1$. This establishes the desired involution. If $i_1\neq 1$ or $j_1\neq 0$, then
$\mathcal{N}_1=\mathcal{N}_2=\emptyset$ and hence $\mathcal{N}(\mathbf{P}_I,\mathbf{P}'_J)=S_{I,J}$; for this subcase, we simply take $\phi$ to be the identity map.

We proceed to define $\phi$ for $k\geq 2$. In this case we divide $\mathcal{N}_1$ into the following two subsets:
\begin{align*}
\mathcal{N}_{1,1}=\{\mathbf{p} \in \mathcal{N}_1 \mid p_2 = P_2 \to Q_1 \to P'_{1}\}, \quad \mathcal{N}_{1,2}= \mathcal{N}_1 \setminus \mathcal{N}_{1,1}.
\end{align*}
In the following, we will define a sign-reversing involution $\phi$ on $\mathcal{N}_1\cup \mathcal{N}_2\cup\mathcal{N}_3$ such that
\begin{align*}
\phi(\mathcal{N}_{1,1})=\mathcal{N}_3, \quad \phi(\mathcal{N}_{1,2})=\mathcal{N}_2.
\end{align*}
There are several subcases to consider.

\begin{itemize}
\item[(i)] If $i_1 = 0$, $i_1 \ge 2$, or $i_1=j_1 = 1$,  then $\mathcal{N}_{1,1} = \mathcal{N}_{1,2} = \mathcal{N}_2 = \mathcal{N}_3 = \emptyset$. For these three situations, take $\phi$ to be the identity map.

\item[(ii)] If $i_1 = 1$, $j_1 = 0$, and moreover $i_2 \ge 3$ or $j_2 \ge 2$, then $\mathcal{N}_{1,1} = \mathcal{N}_3 = \emptyset$ and $\mathcal{N}_{1,2} = \mathcal{N}_1$. For these situations, we take $\phi$
    to be the map which sends
\[
(P_1 \to P'_{0},p_2,\ldots,p_k)\in \mathcal{N}_{1,2} \quad \mbox{and} \quad
(P_1 \to Q_1 \overset{l}{\to} P'_{0},p_2,\ldots,p_k)\in \mathcal{N}_{2}
\]
to each other. It is clear that $(P_1 \to P'_{0},p_2,\ldots,p_k)$ and $(P_1 \to Q_1 \overset{l}{\to} P'_{0},p_2,\ldots,p_k)$  have opposite weights.

\item[(iii)]  If $i_1 = 1$, $j_1 = 0$ and $i_2 = 2$, $j_2 = 1$, then
\begin{align*}
  \mathcal{N}_{1,1}&=\{\mathbf{p}\in \mathcal{N}(P_I,P'_J) \mid p_1=P_1 \to P'_{0}, p_2 = P_2 \to Q_1 \to P'_{1}\}, \\
  \mathcal{N}_{1,2}&=\{\mathbf{p} \in \mathcal{N}(P_I,P'_J) \mid p_1=P_1 \to P'_{0}, p_2 = P_2 \to Q_2 \to P'_{1}\}, \\
  \mathcal{N}_{2}&=\{\mathbf{p} \in \mathcal{N}(P_I,P'_J) \mid p_1 = P_1 \to Q_1 \overset{l}{\to} P'_0, p_2 = P_2 \to Q_2 \to P'_{1}\},\\
  \mathcal{N}_{3}&=\left\{\mathbf{p} \in \mathcal{N}(P_I,P'_J)\quad
  \left | \begin{array}{ll}
  &p_1 = P_1 \to Q_0 \to P'_0,\\
  & \exists\, l \ge 2 \mbox{ such that } p_l = P_l \to Q_l \to P'_{l-1},
  \\
  &p_m = P_m \to Q_{m-1} \to P'_{m-1}, \forall\, 2 \le m \le l-1
  \end{array}\right.
  \right\}.
\end{align*}
Now we are going to define the map $\phi$ on $\mathcal{N}_1\cup \mathcal{N}_2\cup\mathcal{N}_3$.

If $\mathbf{p}\in \mathcal{N}_{1,1}$, say
$\mathbf{p}=(P_1 \to P'_{0},p_2,\ldots,p_k)$, then we take $l$  to be the largest number such that $p_m = P_m \to Q_{m-1} \to P'_{m-1}$ for $2 \le m \le l$, and let
\begin{align*}
\phi(\mathbf{p})&=(P_1 \to Q_0 \to P'_0,p_2,\ldots,p_{l-1},P_l \to Q_l \to P'_{l-1},p_{l+1},\ldots,p_k).
\end{align*}
Thus, $\phi(\mathbf{p})\in \mathcal{N}_3$.

If $\mathbf{p}\in \mathcal{N}_{3}$, namely, there exists $l \ge 2$ such that $\mathbf{p} = (P_1 \to Q_0 \to P'_0,p_2,\ldots,p_k)$ with $p_l = P_l \to Q_l \to P'_{l-1}$ and $p_m = P_m \to Q_{m-1} \to P'_{m-1}$ for any $2 \le m \le l-1 $, then let
\begin{align*}
\phi(\mathbf{p})&=(P_1 \to P'_0,p_2,\ldots,p_{l-1},P_l \to Q_{l-1} \to P'_{l-1},p_{l+1},\ldots,p_k).
\end{align*}
Here the map $\phi$ is well-defined since the number $l$ exists then it must be unique by the definition of $(\mathcal{P}_3)$. It is also clear that $\phi(\mathbf{p})\in \mathcal{N}_{1,1}$.

If $\mathbf{p}\in \mathcal{N}_{1,2}$, say
$\mathbf{p}=(P_1 \to P'_{0},p_2,\ldots,p_k)$, then let
\begin{align*}
\phi(\mathbf{p})&=(P_1 \to Q_1 \overset{l}{\to} P'_0,p_2,\ldots,p_k).
\end{align*}
Hence, we have $\phi(\mathbf{p})\in \mathcal{N}_{2}$.

If $\mathbf{p}\in \mathcal{N}_{2}$, say
$\mathbf{p}=(P_1 \to Q_1 \overset{l}{\to} P'_0,p_2,\ldots,p_k)$, then let
\begin{align*}
\phi(\mathbf{p})&=(P_1 \to P'_{0},p_2,\ldots,p_k).
\end{align*}
It is obvious that $\phi(\mathbf{p})\in \mathcal{N}_{1,2}$.

Figure \ref{fig-phi-iii} gives an illustration of $\phi$ for this subcase.

\begin{figure}[ht]
  \centering
\begin{tikzpicture}
[place/.style={thick,fill=black!100,circle,inner sep=0pt,minimum size=1mm,draw=black!100},scale=0.9]
\node [place,label=below:{\footnotesize$P_1$}] (v1) at (-2,0.5) {};
\node [place,label=below:{\footnotesize$P'_0$}] (v3) at (0,1) {};
\node [place,label=below:{\footnotesize$P_2$}] (v4) at (-2,-0.5) {};
\node [place,label=above:{\footnotesize$Q_1$}] (v5) at (-1,0) {};
\node [place,label=below:{\footnotesize$P'_1$}] (v6) at (0,0) {};
\node at (-1,-0.5) {$\vdots$};
\node [place,label=below:{\footnotesize$P_{l-1}$}] (v7) at (-2,-2) {};
\node [place,label=above:{\footnotesize$Q_{l-2}$}] (v8) at (-1,-1.5) {};
\node [place,label=below:{\footnotesize$P'_{l-2}$}] (v9) at (0,-1.5) {};
\node [place,label=below:{\footnotesize\footnotesize{$P_l$}}] (v12) at (-2,-3) {};
\node [place,label=above:{\footnotesize$Q_{l-1}$}] (v10) at (-1,-2.5) {};
\node [place,label=below:{\footnotesize$P'_{l-1}$}] (v11) at (0,-2.5) {};
\node at (-1,-4) {$\mathcal{N}_{1,1}$};
\draw [thick] (v1) -- (v3);
\draw [thick] (v4) -- (v5) -- (v6);
\draw [thick] (v7) -- (v8) -- (v9);
\draw [thick] (v12) -- (v10) -- (v11);
\node at (1.5,-0.5) {$\overset{\phi}{\leftrightarrow}$};
\node [place,label=below:{\footnotesize$P_1$}] (v13) at (3,0.5) {};
\node [place,label=above:{\footnotesize$Q_0$}] (v2) at (4,1) {};
\node [place,label=below:{\footnotesize$P'_0$}] (v14) at (5,1) {};
\node [place,label=below:{\footnotesize$P_2$}] (v15) at (3,-0.5) {};
\node [place,label=above:{\footnotesize$Q_1$}] (v16) at (4,0) {};
\node [place,label=below:{\footnotesize$P'_1$}] (v17) at (5,0) {};
\node at (4,-0.5) {$\vdots$};
\node [place,label=below:{\footnotesize$P_{l-1}$}] (v18) at (3,-2) {};
\node [place,label=above:{\footnotesize$Q_{l-2}$}] (v19) at (4,-1.5) {};
\node [place,label=below:{\footnotesize$P'_{l-2}$}] (v20) at (5,-1.5) {};
\node [place,label=below:{\footnotesize$P_l$}] (v21) at (3,-3) {};
\node [place,label=below:{\footnotesize$Q_{l}$}] (v22) at (4,-3) {};
\node [place,label=below:{\footnotesize$P'_{l-1}$}] (v23) at (5,-2.5) {};
\node at (4,-4) {$\mathcal{N}_{3}$};
\draw [thick] (v13) -- (v2) -- (v14);
\draw [thick] (v15) -- (v16) -- (v17);
\draw [thick] (v18) -- (v19) -- (v20);
\draw [thick] (v21) -- (v22) -- (v23);
\node [place,label=below:{\footnotesize$P_1$}] (v24) at (7.5,0.5) {};
\node [place,label=below:{\footnotesize$P'_0$}] (v25) at (9.5,1) {};
\node [place,label=below:{\footnotesize$P_2$}] (v26) at (7.5,-0.5) {};
\node[place,label=below:{\footnotesize$Q_2$}] (v27) at (8.5,-0.5) {};
\node[place,label=below:{\footnotesize$P'_1$}] (v28) at (9.5,0) {};
\node at (7.5,-1.75) {$\vdots$};
\node at (8.5,-1.75) {$\vdots$};
\node at (9.5,-1.75) {$\vdots$};
\node at (8.5,-4) {$\mathcal{N}_{1,2}$};
\draw [thick](v24) -- (v25);
\draw[thick] (v26) -- (v27) -- (v28);
\node at (10.75,-0.5) {$\overset{\phi}{\leftrightarrow}$};
\node [place,label=below:{\footnotesize$P_1$}] (v35) at (12,0.5) {};
\node[place,label=below:{\footnotesize$Q_1$}]  (v36) at (13,0.5) {};
\node [place,label=below:{\footnotesize$P'_0$}]at (14,1) {};
\node[place,label=below:{\footnotesize$P_2$}] (v37) at (12,-0.5) {};
\node[place,label=below:{\footnotesize$Q_2$}]  (v38) at (13,-0.5) {};
\node[place,label=below:{\footnotesize$P'_1$}]  (v39) at (14,0) {};
\node at (12,-1.75) {$\vdots$};
\node at (13,-1.75) {$\vdots$};
\node at (14,-1.75) {$\vdots$};
\node at (13,-4) {$\mathcal{N}_{2}$};
\draw[thick] (v35) -- (v36);
\draw[thick] (v37) -- (v38) -- (v39);
\draw [thick]plot[smooth, tension=.7] coordinates {(13,0.5) (13.4,0.85) (14,1)};
\draw [thick][->](13.4,0.85)-- (13.41,0.86);
\end{tikzpicture}\caption{An illustration of $\phi$ for subcase (iii)}\label{fig-phi-iii}
\end{figure}
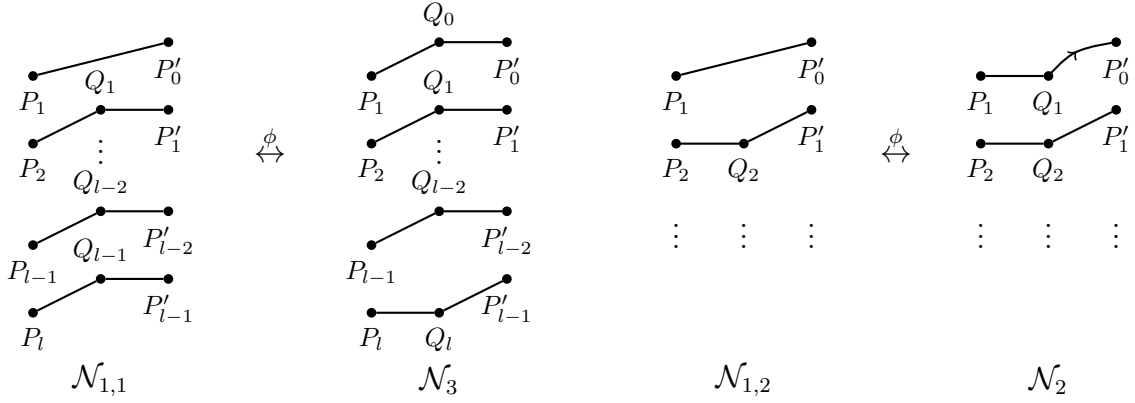

\end{itemize}

With the above definition of $\phi$, it is straightforward to verify that
$\phi$ is a sign-reversing involution on $\mathcal{N}_{1}\cup\mathcal{N}_{2}\cup\mathcal{N}_3$, as desired. \qed

\section{The Hankel matrix \texorpdfstring{$H$}{}}\label{sect-pf}
The aim of this section is to give a combinatorial proof of the $q$-total positivity of the Hankel matrix $H$ of Narayana polynomials of type $B$.
It suffices to combinatorially prove the $q$-total positivity of each
leading principal submatrix of $H$. To this end, let us first establish the planar networks for its leading principal submatrices.

Let $L^B_n=(l_{i,j}(q))_{0\le i,j \le n+1}$, $B_n=(b_{i,j}(q))_{0\le i,j \le n}$ and $H_n = (b_{i+j,0}(q))_{0 \le i,j \le n}$. By \eqref{eq-narab-recurrence} it is evident that
\begin{align}\label{eq-recurrence}
B_{n+1}=\bar{B}_nL^B_{n}, \mbox{ where } \bar{B}_n=
\begin{pmatrix}
1 & O\\
O & B_n
\end{pmatrix}.
\end{align}
Aigner \cite{Aig01} proved that
\begin{equation}\label{eq-Hankel-recurrence}
H_n=B_nT_nB_n^{t},
\end{equation}
where
\begin{equation*}
T_n = \begin{pmatrix}
1& & & & \\
 &2q& & & \\
 & &2q^2& & \\
 & & &\ddots& \\
 & & & &2q^n
\end{pmatrix}
_{(n+1) \times (n+1)},
\end{equation*}
and $B_n^t$ denotes the transpose of $B_n$.

These two formulas allow us to recursively construct the planar network for $H_n$. Precisely, we mainly make use of the following lemma, which provides a way to build a network for a product of matrices. Recall that in a digraph, a vertex $v$ is called a source (resp. sink) if there is no arcs point in (resp. out of) it. The following result could be considered as a corollary of the transfer-matrix method (see \cite[Theorem
4.7.1]{StaEC1} for instance), but for self-containedness we will give a detailed proof.

\begin{lem}\label{lem-transfer-matrix}
Given $k$ square matrices $M_1$, $M_2$, \ldots, $M_k$ of order $n$, for each $1 \le i \le k$ assume that $\mathcal{M}_i = (D^{M_i},\mathrm{wt}_{D^{M_i}},(u_{i,1},\ldots,u_{i,n}),
(v_{i,1},\ldots,v_{i,n}))$ is a planar network for $M_i$ with $u_{i,j}$ being a source and $v_{i,j}$ being a sink for all $1 \le j \le n$. Let $D^M$ be the digraph obtained by placing $D^{M_1}$, $D^{M_2}$, \ldots, $D^{M_k}$ in succession and identifying $v_{i,j}$ with $u_{i+1,j}$ for each $1 \le i \le k-1$ and $1 \le j \le n$, and let $\mathrm{wt}_{D^M}$ be the weight function inherited from $\mathrm{wt}_{D^{M_1}},\ldots,\mathrm{wt}_{D^{M_k}}$ in an obvious way. Then
\[
\mathcal{M} = (D^M, \mathrm{wt}_{D^M}, (u_{1,1},\ldots,u_{1,n}), (v_{k,1},\ldots,v_{k,n}))
\]
is a planar network for the product $M_1\cdots M_k$.
\end{lem}

\pf The proof is by induction on $k$. Let us first prove the base case $k = 2$. The construction tells that for any $1 \le i,j \le n$, each directed path from $u_{1,i}$ to $v_{2,j}$ must pass through exactly one vertex $v_{1,l}$ ($=u_{2,l}$) for some $1 \le l \le n$, and hence
\begin{equation*}
GF_{D^M}(u_{1,i},v_{2,j}) = \sum_{l=1}^{n} GF_{D^{M_1}}(u_{1,i},v_{1,l}) GF_{D^{M_2}}(u_{2,l},v_{2,j}),
\end{equation*}
as desired. Assume the assertion for $k$ ($k \ge 2)$. By applying the preceding proof to $M_1\cdots M_k$ and $M_{k+1}$, we find that the assertion also holds for $k+1$. This completes the proof. \qed

\begin{figure}[htb]
  \centering
\begin{tikzpicture}
[place/.style={thick,fill=black!100,circle,inner sep=0pt,minimum size=1mm,draw=black!100},scale=1.5]
\draw [thick] [->] (2.5,2.5) -- (3,2.5);
\draw [thick] [->] (3,2.5) -- (3.5,2.5) -- (4,2.5);
\draw [thick] (4,2.5) -- (4.5,2.5);
\node [place,label=above:{\footnotesize$P_3^{(2)}$}] at (2.5,2.5) {};
\node [place,label=above:{\footnotesize$Q_3^{(2)}$}] at (3.5,2.5) {};
\node [place,label=above:{\footnotesize$P_3^{(3)}$}] at (4.5,2.5) {};
\draw [thick] [->](2.5,1.5) -- (3,1.5);
\draw [thick] [->] (3,1.5) -- (4,1.5);
\draw [thick] (4,1.5) -- (4.5,1.5);
\node [place,label=below:{\footnotesize$P_2^{(2)}$}] at (2.5,1.5) {};
\node [place,label=below:{\footnotesize$Q_2^{(2)}$}] at (3.5,1.5) {};
\node [place,label=below:{\footnotesize$P_2^{(3)}$}] at (4.5,1.5) {};
\draw [thick] [->] (2.5,1.5) --(3,2);
\draw [thick] (3,2) --  (3.5,2.5) ;
\draw [thick] [->] (2.5,1.5) -- (3.5,2);
\draw [thick] (3.5,2) -- (4.5,2.5);
\draw [thick]plot[smooth, tension=.7] coordinates {(3.5,1.5) (4.1,1.9) (4.5,2.5)};
\draw [thick]plot[smooth, tension=.7] coordinates {(3.5,1.5) (3.9,2.1) (4.5,2.5)};
\draw [thick][->] (4.1,1.9) -- (4.11,1.91);
\draw [thick][->] (3.9,2.1) -- (3.91,2.11);
\draw [thick] [->] (2.5,0.5) -- (3,0.5);
\draw [thick] [->] (3,0.5) -- (4,0.5);
\draw [thick] (4,0.5) -- (4.5,0.5);
\node [place,label=below:{\footnotesize$P_1^{(2)}$}] at (2.5,0.5) {};
\node [place,label=below:{\footnotesize$Q_1^{(2)}$}] at (3.5,0.5) {};
\node [place,label=below:{\footnotesize$P_1^{(3)}$}] at (4.5,0.5) {};
\draw [thick] [->] (2.5,0.5) -- (3,1);
\draw [thick] (3,1) -- (3.5,1.5);
\draw [thick] [->] (3.5,0.5) -- (4,1);
\draw [thick] (4,1) -- (4.5,1.5);
\draw [thick] [->] (2.5,-0.5) -- (3,-0.5);
\draw [thick] [->] (3,-0.5) -- (4,-0.5);
\draw [thick]  (4,-0.5) -- (4.5,-0.5);
\node [place,label=below:{\footnotesize$P_0^{(2)}$}] at (2.5,-0.5) {};
\node [place,label=below:{\footnotesize$Q_0^{(2)}$}] at (3.5,-0.5) {};
\node [place,label=below:{\footnotesize$P_0^{(3)}$}] at (4.5,-0.5) {};
\draw [thick] [->] (2.5,-0.5) -- (3,0);
\draw [thick] (3,0) -- (3.5,0.5);
\draw [thick] [->] (3.5,-0.5) -- (4,0);
\draw [thick] (4,0) -- (4.5,0.5);
\node [blue] at (2.75,2) {$q$};
\node [blue] at (3.375,1.75) {$-1$};
\node [blue] at (2.75,1) {$q$};
\node [blue] at (2.75,0) {$q$};
\end{tikzpicture}
\caption{Planar network $\mathcal{L}^{\mathcal{B}}_2$}\label{fig-network-lb2}
\end{figure}
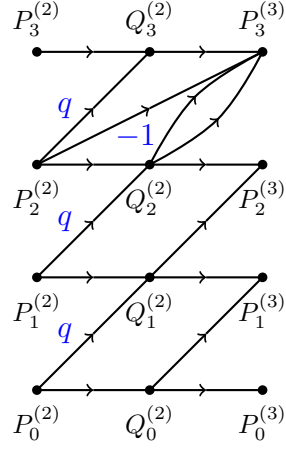

Now we present the construction for the planar network $\mathcal{H}_n$ for $H_n$, which is essentially based on $\mathcal{L^B}$. At first, we give the planar network for $L^B_n$, which is actually obtained by cutting off the part of $\mathcal{L^B}$ below $y = -n-1$. Precisely, let $D^{L^B_n}$ be the subgraph of $D^{L^B}$ induced by the vertices $P_0,\ldots,P_{n+1},Q_0,\ldots,Q_{n+1},P'_0,\ldots,P'_{n+1}$, and let $\mathrm{wt}_{D^{L^B_n}}$ be the restriction of $\mathrm{wt}_{D^{L^B}}$ to $D^{L^B_n}$. Then
$(D^{L^B_n},\mathrm{wt}_{D^{L^B_n}},(P_0,\ldots,P_{n+1}),(P'_0,\ldots,P'_{n+1}))$ is a planar network for $L^B_n$. Unfortunately, this labeling is not convenient for introducing the recursive construction of $\mathcal{H}_n$. In the rest of this paper, we will label the vertex $P_i$ by $P^{(n)}_{n+1-i}$, $Q_i$ by $Q^{(n)}_{n+1-i}$, $P'_i$ by $P^{(n+1)}_{n+1-i}$ for $0 \le i \le n+1$ in the digraph $D^{L^B_n}$. Moreover, we may shift the digraphs $D^{L^B_n}$ in the plane such that $P^{(i)}_j = (2i,j)$ and $Q^{(i)}_j = (2i+1,j)$ for all $i,j \ge 0$.
Then
\[
\mathcal{L}^{\mathcal{B}}_n = (D^{L^B_n},\mathrm{wt}_{D^{L^B_n}},(P^{(n)}_{n+1},\ldots,P^{(n)}_0),(P^{(n+1)}_{n+1},\ldots,P^{(n+1)}_0))
\]
is a planar network for $L^B_n$, and
\[
L^B_n = \left(GF_{D^{L^B_n}}(P^{(n)}_{n+1-i},P^{(n+1)}_{n+1-j})\right)_{0 \le i,j \le n+1}.
\]
Figure \ref{fig-network-lb2} shows the planar network $\mathcal{L}^{\mathcal{B}}_2$.

Using the planar networks $\mathcal{L}^{\mathcal{B}}_0,\mathcal{L}^{\mathcal{B}}_1,\ldots,\mathcal{L}^{\mathcal{B}}_{n}$, we can construct the planar network $\mathcal{B}_{n+1}$ for $B_{n+1}$.

\begin{itemize}
\item For $n=0$, we take $\mathcal{B}_{1}$ to be the planar network $\mathcal{L}^{\mathcal{B}}_{0}$ since $B_1 = L^B_0$.


\item Assuming that $\mathcal{B}_{n}$
has been constructed for some $n \geq 1$, we continue to construct $\mathcal{B}_{n+1}$. Let $D^{\bar{B}_n}$ be the digraph with $V(D^{\bar{B}_n}) = V(D^{B_n}) \cup \{P_{n+1}^{(0)},P_{n+1}^{(1)},\ldots,P_{n+1}^{(n)}\}$ and $A(D^{\bar{B}_n}) = A(D^{B_n}) \cup \{P_{n+1}^{(i)} \to P_{n+1}^{(i+1)} \mid 0\leq i\leq n-1\}$, and let $\mathrm{wt}_{D^{\bar{B}_n}}(a)$ be equal to $\mathrm{wt}_{D^{B_n}}(a)$ for $a \in A(D^{B_n})$ and equal to 1 for the other arcs. Then
\[
\mathcal{\bar{B}}_n = (D^{\bar{B}_n},\mathrm{wt}_{D^{\bar{B}_n}}, (P_{n+1}^{(0)},P_n^{(0)},\ldots,P_0^{(0)}), (P_{n+1}^{(n)},P_n^{(n)},\ldots,P_0^{(n)}))
\]
is a planar network for $\bar{B}_n$.
By \eqref{eq-recurrence} and Lemma \ref{lem-transfer-matrix}, we obtain that
\[
\mathcal{B}_{n+1} = (D^{B_{n+1}},\mathrm{wt}_{D^{B_{n+1}}}, (P_{n+1}^{(0)},P_{n}^{(0)},\ldots,P_{0}^{(0)}), (P_{n+1}^{(n+1)},P_{n}^{(n+1)},\ldots,P_{0}^{(n+1)}))
\]
is a planar network for $B_{n+1}$, where $D^{B_{n+1}}$ and $\mathrm{wt}_{D^{B_{n+1}}}$ are defined in the way as described in Lemma \ref{lem-transfer-matrix}. See Figure \ref{fig-db3} for an illustration of $D^{B_3}$.

\begin{figure}[htb]
\centering
\begin{tikzpicture}
[place/.style={thick,fill=black!100,circle,inner sep=0pt,minimum size=1mm,draw=black!100},scale=1.5]
\draw [thick] [->] (-1.5,2.5) -- (-0.5,2.5);
\draw [thick] [->] (-0.5,2.5) -- (1.5,2.5);
\draw [thick] (1.5,2.5) -- (2.5,2.5);
\draw [thick] [->] (2.5,2.5) -- (3,2.5);
\draw [thick] [->] (3,2.5) -- (3.5,2.5) -- (4,2.5);
\draw [thick] (4,2.5) -- (4.5,2.5);
\node [place,label=below:{\footnotesize$P_3^{(0)}$}] at (-1.5,2.5) {};
\node [place,label=below:{\footnotesize$P_3^{(1)}$}] at (0.5,2.5) {};
\node [place,label=below:{\footnotesize$P_3^{(2)}$}] at (2.5,2.5) {};
\node [place,label=below:{\footnotesize$Q_3^{(2)}$}] at (3.5,2.5) {};
\node [place,label=below:{\footnotesize$P_3^{(3)}$}] at (4.5,2.5) {};
\draw [thick] [->] (-1.5,1.5) -- (-0.5,1.5);
\draw [thick] [->] (-0.5,1.5) -- (1,1.5);
\draw [thick] [->] (1,1.5) -- (2,1.5);
\draw [thick] (2,1.5) --(2.5,1.5);
\draw [thick][->] (2.5,1.5) -- (3,1.5);
\draw [thick] (3,1.5) --(3.5,1.5);
\draw [thick][->] (3.5,1.5) -- (4,1.5);
\draw [thick] (4,1.5) -- (4.5,1.5);
\node [place,label=below:{\footnotesize$P_2^{(0)}$}] at (-1.5,1.5) {};
\node [place,label=below:{\footnotesize$P_2^{(1)}$}] at (0.5,1.5) {};
\node [place,label=below:{\footnotesize$Q_2^{(1)}$}] at (1.5,1.5) {};
\node [place,label=below:{\footnotesize$P_2^{(2)}$}] at (2.5,1.5) {};
\node [place,label=below:{\footnotesize$Q_2^{(2)}$}] at (3.5,1.5) {};
\node [place,label=below:{\footnotesize$P_2^{(3)}$}] at (4.5,1.5) {};
\draw [thick] [->] (2.5,1.5) --(3,2);
\draw [thick]  (3,2) --  (3.5,2.5) ;
\draw [thick] [->] (2.5,1.5) -- (3.5,2);
\draw [thick]  (3.5,2) -- (4.5,2.5);
\draw [thick] [->] (-1.5,0.5) -- (-1,0.5);
\draw [thick] [->] (-1,0.5) --  (0,0.5);
\draw [thick] (0,0.5) -- (0.5,0.5);
\draw [thick] [->] (0.5,0.5) -- (1,0.5);
\draw [thick] [->] (1,0.5) -- (2,0.5);
\draw [thick]  (2,0.5) -- (2.5,0.5);
\draw [thick] [->] (2.5,0.5) -- (3,0.5);
\draw [thick] [->] (3,0.5) -- (4,0.5);
\draw [thick] (4,0.5) -- (4.5,0.5);
\node [place,label=below:{\footnotesize$P_1^{(0)}$}] at (-1.5,0.5) {};
\node [place,label=below:{\footnotesize$Q_1^{(0)}$}] at (-0.5,0.5) {};
\node [place,label=below:{\footnotesize$P_1^{(1)}$}] at (0.5,0.5) {};
\node [place,label=below:{\footnotesize$Q_1^{(1)}$}] at (1.5,0.5) {};
\node [place,label=below:{\footnotesize$P_1^{(2)}$}] at (2.5,0.5) {};
\node [place,label=below:{\footnotesize$Q_1^{(2)}$}] at (3.5,0.5) {};
\node [place,label=below:{\footnotesize$P_1^{(3)}$}] at (4.5,0.5) {};
\draw [thick] [->] (0.5,0.5) --(1,1);
\draw [thick] (1,1) -- (1.5,1.5);
\draw [thick] (1.5,1) -- (2.5,1.5);
\draw [thick] [->]  (2.5,0.5) -- (3,1);
\draw [thick]  (3,1) -- (3.5,1.5);
\draw [thick] [->] (3.5,0.5) -- (4,1);
\draw [thick] (4,1) -- (4.5,1.5);
\draw [thick] [->] (-1.5,-0.5) -- (-1,-0.5);
\draw [thick] [->] (-1,-0.5) -- (-0.5,-0.5) -- (0,-0.5);
\draw [thick] [->] (0,-0.5) -- (1,-0.5);
\draw [thick] [->] (1,-0.5) -- (1.5,-0.5) -- (2,-0.5);
\draw [thick] [->] (2,-0.5) -- (2.5,-0.5) -- (3,-0.5);
\draw [thick] [->] (3,-0.5) -- (3.5,-0.5) -- (4,-0.5);
\draw [thick] (4,-0.5) -- (4.5,-0.5);
\node [place,label=below:{\footnotesize$P_0^{(0)}$}] at (-1.5,-0.5) {};
\node [place,label=below:{\footnotesize$Q_0^{(0)}$}] at (-0.5,-0.5) {};
\node [place,label=below:{\footnotesize$P_0^{(1)}$}] at (0.5,-0.5) {};
\node [place,label=below:{\footnotesize$Q_0^{(1)}$}] at (1.5,-0.5) {};
\node [place,label=below:{\footnotesize$P_0^{(2)}$}] at (2.5,-0.5) {};
\node [place,label=below:{\footnotesize$Q_0^{(2)}$}] at (3.5,-0.5) {};
\node [place,label=below:{\footnotesize$P_0^{(3)}$}] at (4.5,-0.5) {};
\draw [thick] [->] (0.5,-0.5) -- (1,0);
\draw [thick] (1,0) --(1.5,0.5);
\draw [thick] [->]  (1.5,-0.5) -- (2,0);
\draw [thick] (2,0) -- (2.5,0.5);
\draw [thick] [->]  (2.5,-0.5) -- (3,0);
\draw [thick] (3,0) -- (3.5,0.5);
\draw [thick] [->] (3.5,-0.5) -- (4,0);
\draw [thick] (4,0) -- (4.5,0.5);
\draw [thick] [->] (-1.5,-0.5)  -- (-1,0);
\draw [thick] (-1,0) -- (-0.5,0.5);
\draw [thick] [->] (-1.5,-0.5) -- (-0.5,0);
\draw [thick]  (-0.5,0) -- (0.5,0.5);
\draw [thick][->] (0.5,0.5) -- (1.5,1);
\draw [thick]plot[smooth, tension=.7] coordinates {(-0.5,-0.5) (0.1,-0.1) (0.5,0.5)};
\draw [thick]plot[smooth, tension=.7] coordinates {(-0.5,-0.5) (-0.1,0.1) (0.5,0.5)};
\draw [thick]plot[smooth, tension=.7] coordinates {(1.5,0.5) (2.1,0.9) (2.5,1.5)};
\draw [thick]plot[smooth, tension=.7] coordinates {(1.5,0.5) (1.9,1.1) (2.5,1.5)};
\draw [thick]plot[smooth, tension=.7] coordinates {(3.5,1.5) (4.1,1.9) (4.5,2.5)};
\draw [thick]plot[smooth, tension=.7] coordinates {(3.5,1.5) (3.9,2.1) (4.5,2.5)};
\draw [thick][->] (0.1,-0.1) -- (0.11,-0.09);
\draw [thick][->] (-0.1,0.1) -- (-0.09,0.11);
\draw [thick][->] (2.1,0.9) -- (2.11,0.91);
\draw [thick][->] (1.9,1.1) -- (1.91,1.11);
\draw [thick][->] (4.1,1.9) -- (4.11,1.91);
\draw [thick][->] (3.9,2.1) -- (3.91,2.11);
\end{tikzpicture}
\caption{Digraph $D^{B_3}$}\label{fig-db3}
\end{figure}
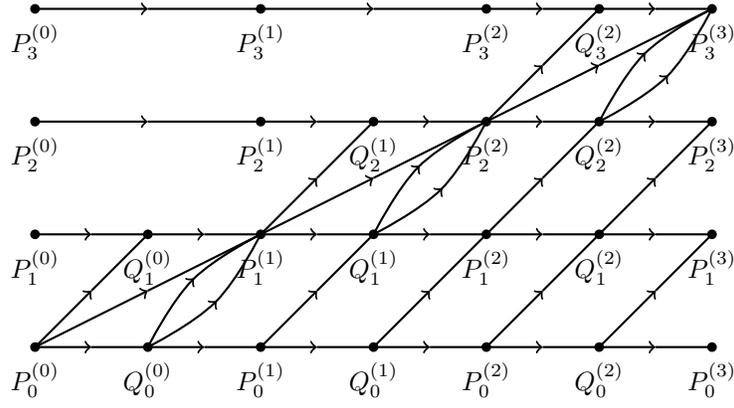
\end{itemize}

\begin{figure}[htb]
\begin{tikzpicture}
[place/.style={thick,fill=black!100,circle,inner sep=0pt,minimum size=1mm,draw=black!100},scale=1.15]
\draw [thick] [->] (-1.5,2.5) -- (-0.5,2.5);
\draw [thick]  (-0.5,2.5) -- (0.5,2.5);
\draw [thick] [->] (0.5,2.5) -- (1.5,2.5);
\draw [thick]  (1.5,2.5) -- (2.5,2.5);
\draw [thick] [->]  (2.5,2.5) -- (3,2.5);
\draw [thick] [->] (3,2.5) -- (3.5,2.5) -- (4,2.5);
\draw [thick] [->] (4,2.5) -- (4.5,2.5) -- (5,2.5);
\draw [thick] [->] (5,2.5) -- (5.5,2.5) -- (6,2.5);
\draw [thick] [->] (6,2.5) -- (6.5,2.5) -- (7,2.5);
\draw [thick] [->]  (7,2.5) -- (8.5,2.5);
\draw [thick] (8.5,2.5) -- (9.5,2.5);
\draw [thick] [->] (9.5,2.5) -- (10.5,2.5);
\draw [thick] (10.5,2.5) -- (11.5,2.5);
\node [place,label=below:{\tiny$P_3^{(0)}$}] at (-1.5,2.5) {};
\node [place,label=below:{\tiny$P_3^{(1)}$}] at (0.5,2.5) {};
\node [place,label=below:{\tiny$P_3^{(2)}$}] at (2.5,2.5) {};
\node [place,label=below:{\tiny$Q_3^{(2)}$}] at (3.5,2.5) {};
\node [place,label=below:{\tiny$P_3^{(3)}$}] at (4.5,2.5) {};
\node [place,label=below:{\tiny$\bar{P}_3^{(3)}$}] at (5.5,2.5) {};
\node [place,label=below:{\tiny$\bar{Q}_3^{(2)}$}] at (6.5,2.55) {};
\node [place,label=below:{\tiny$\bar{P}_3^{(2)}$}] at (7.5,2.5) {};
\node [place,label=below:{\tiny$\bar{P}_3^{(1)}$}] at (9.5,2.5) {};
\node [place,label=below:{\tiny$\bar{P}_3^{(0)}$}] at (11.5,2.5) {};
\draw [thick]plot[smooth, tension=.7] coordinates {(5.5,2.5) (6.1,2.1) (6.5,1.5)};
\draw [thick]plot[smooth, tension=.7] coordinates {(5.5,2.5) (5.9,1.9) (6.5,1.5)};
\draw [thick] [->] (6.09,2.11) -- (6.1,2.1);
\draw [thick] [->] (5.89,1.91) -- (5.9,1.9);
\draw [thick] [->] (5.5,2.5) -- (6.5,2);
\draw [thick] (6.5,2) -- (7.5,1.5);
\draw [thick] [->] (6.5,2.5) -- (7,2);
\draw [thick] (7,2) -- (7.5,1.5);
\draw [thick] [->] (-1.5,1.5) -- (-0.5,1.5);
\draw [thick]  (-0.5,1.5) -- (0.5,1.5);
\draw [thick] [->]  (0.5,1.5) -- (1,1.5);
\draw [thick] [->] (1,1.5) -- (1.5,1.5) -- (2,1.5);
\draw [thick] [->] (2,1.5) -- (2.5,1.5) -- (3,1.5);
\draw [thick] [->] (3,1.5) -- (3.5,1.5) -- (4,1.5);
\draw [thick] [->] (4,1.5) -- (4.5,1.5) -- (5,1.5);
\draw [thick] [->] (5,1.5) -- (5.5,1.5) -- (6,1.5);
\draw [thick] [->] (6,1.5) -- (6.5,1.5) -- (7,1.5);
\draw [thick] [->] (7,1.5) -- (7.5,1.5) -- (8,1.5);
\draw [thick] [->] (8,1.5) -- (8.5,1.5) -- (9,1.5);
\draw [thick] [->] (9,1.5) --  (10.5,1.5);
\draw [thick] (10.5,1.5) --  (11.5,1.5);
\node [place,label=below:{\tiny$P_2^{(0)}$}] at (-1.5,1.5) {};
\node [place,label=below:{\tiny$P_2^{(1)}$}] at (0.5,1.5) {};
\node [place,label=below:{\tiny$Q_2^{(1)}$}] at (1.5,1.5) {};
\node [place,label=below:{\tiny$P_2^{(2)}$}] at (2.5,1.5) {};
\node [place,label=below:{\tiny$Q_2^{(2)}$}] at (3.5,1.5) {};
\node [place,label=below:{\tiny$P_2^{(3)}$}] at (4.5,1.5) {};
\node [place,label=below:{\tiny$\bar{P}_2^{(3)}$}] at (5.5,1.5) {};
\node [place,label=below:{\tiny$\bar{Q}_2^{(2)}$}] at (6.5,1.5) {};
\node [place,label=below:{\tiny$\bar{P}_2^{(2)}$}] at (7.5,1.5) {};
\node [place,label=below:{\tiny$\bar{Q}_2^{(1)}$}] at (8.5,1.55) {};
\node [place,label=below:{\tiny$\bar{P}_2^{(1)}$}] at (9.5,1.5) {};
\node [place,label=below:{\tiny$\bar{P}_2^{(0)}$}] at (11.5,1.5) {};
\draw [thick] [->] (2.5,1.5) --(3,2);
\draw [thick]  (3,2) --  (3.5,2.5) ;
\draw [thick] [->] (2.5,1.5) -- (3.5,2);
\draw [thick] (3.5,2) -- (4.5,2.5);
\draw [thick]plot[smooth, tension=.7] coordinates {(3.5,1.5) (4.1,1.9) (4.5,2.5)};
\draw [thick]plot[smooth, tension=.7] coordinates {(3.5,1.5) (3.9,2.1) (4.5,2.5)};
\draw [thick][->] (4.1,1.9) -- (4.11,1.91);
\draw [thick][->] (3.9,2.1) -- (3.91,2.11);
\draw [thick] [->] (5.5,1.5) --(6,1);
\draw [thick] (6,1) -- (6.5,0.5);
\draw [thick] [->] (6.5,1.5) -- (7,1);
\draw [thick] (7,1) -- (7.5,0.5);
\draw [thick]plot[smooth, tension=.7] coordinates {(7.5,1.5) (8.1,1.1) (8.5,0.5)};
\draw [thick]plot[smooth, tension=.7] coordinates {(7.5,1.5) (7.9,0.9) (8.5,0.5)};
\draw [thick] [->] (8.09,1.11) -- (8.1,1.1);
\draw [thick] [->] (7.89,0.91) -- (7.9,0.9);
\draw [thick] [->] (7.5,1.5) -- (8.5,1);
\draw [thick] (8.5,1) -- (9.5,0.5);
\draw [thick] [->] (8.5,1.5) -- (9,1);
\draw [thick] (9,1) -- (9.5,0.5);
\draw [thick] [->] (-1.5,0.5) -- (-1,0.5);
\draw [thick] [->] (-1,0.5) -- (-0.5,0.5) -- (0,0.5);
\draw [thick] [->] (0,0.5) -- (0.5,0.5) -- (1,0.5);
\draw [thick] [->] (1,0.5) -- (1.5,0.5) -- (2,0.5);
\draw [thick] [->] (2,0.5) -- (2.5,0.5) -- (3,0.5);
\draw [thick] [->] (3,0.5) -- (3.5,0.5) -- (4,0.5);
\draw [thick] [->] (4,0.5) -- (4.5,0.5) -- (5,0.5);
\draw [thick] [->] (5,0.5) -- (5.5,0.5) -- (6,0.5);
\draw [thick] [->] (6,0.5) -- (6.5,0.5) -- (7,0.5);
\draw [thick] [->] (7,0.5) -- (7.5,0.5) -- (8,0.5);
\draw [thick] [->] (8,0.5) -- (8.5,0.5) -- (9,0.5);
\draw [thick] [->] (9,0.5) -- (9.5,0.5) -- (10,0.5);
\draw [thick] (10,0.5) -- (10.5,0.5);
\draw [thick] [->] (10.5,0.5) -- (11,0.5);
\draw [thick] (11,0.5) -- (11.5,0.5);
\node [place,label=below:{\tiny$P_1^{(0)}$}] at (-1.5,0.5) {};
\node [place,label=below:{\tiny$Q_1^{(0)}$}] at (-0.5,0.5) {};
\node [place,label=below:{\tiny$P_1^{(1)}$}] at (0.5,0.5) {};
\node [place,label=below:{\tiny$Q_1^{(1)}$}] at (1.5,0.5) {};
\node [place,label=below:{\tiny$P_1^{(2)}$}] at (2.5,0.5) {};
\node [place,label=below:{\tiny$Q_1^{(2)}$}] at (3.5,0.5) {};
\node [place,label=below:{\tiny$P_1^{(3)}$}] at (4.5,0.5) {};
\node [place,label=below:{\tiny$\bar{P}_1^{(3)}$}] at (5.5,0.5) {};
\node [place,label=below:{\tiny$\bar{Q}_1^{(2)}$}] at (6.5,0.5) {};
\node [place,label=below:{\tiny$\bar{P}_1^{(2)}$}] at (7.5,0.5) {};
\node [place,label=below:{\tiny$\bar{Q}_1^{(1)}$}] at (8.5,0.5) {};
\node [place,label=below:{\tiny$\bar{P}_1^{(1)}$}] at (9.5,0.5) {};
\node [place,label=below:{\tiny$\bar{Q}_1^{(0)}$}] at (10.5,0.55) {};
\node [place,label=below:{\tiny$\bar{P}_1^{(0)}$}] at (11.5,0.5) {};
\draw [thick] [->] (0.5,0.5) --(1,1);
\draw [thick] (1,1) -- (1.5,1.5);
\draw [thick] [->] (0.5,0.5) -- (1.5,1);
\draw [thick] (1.5,1) -- (2.5,1.5);
\draw [thick]plot[smooth, tension=.7] coordinates {(1.5,0.5) (2.1,0.9) (2.5,1.5)};
\draw [thick]plot[smooth, tension=.7] coordinates {(1.5,0.5) (1.9,1.1) (2.5,1.5)};
\draw [thick][->] (2.1,0.9) -- (2.11,0.91);
\draw [thick][->] (1.9,1.1) -- (1.91,1.11);
\draw [thick] [->] (2.5,0.5) -- (3,1);
\draw [thick] (3,1) -- (3.5,1.5);
\draw [thick] [->] (3.5,0.5) -- (4,1);
\draw [thick] (4,1) -- (4.5,1.5);
\draw [thick] [->] (-1.5,-0.5) -- (-1,-0.5);
\draw [thick] [->] (-1,-0.5) --  (0,-0.5);
\draw [thick] [->] (0,-0.5) -- (1,-0.5);
\draw [thick] [->] (1,-0.5) -- (2,-0.5);
\draw [thick] [->] (2,-0.5) --  (3,-0.5);
\draw [thick] [->] (3,-0.5) --  (4,-0.5);
\draw [thick] [->] (4,-0.5) --  (5,-0.5);
\draw [thick] [->] (5,-0.5) --  (6,-0.5);
\draw [thick] [->] (6,-0.5) -- (7,-0.5);
\draw [thick] [->] (7,-0.5) --  (8,-0.5);
\draw [thick] [->] (8,-0.5) --  (9,-0.5);
\draw [thick] [->] (9,-0.5) -- (10,-0.5);
\draw [thick] (10,-0.5) -- (10.5,-0.5);
\draw [thick] [->] (10.5,-0.5) -- (11,-0.5);
\draw [thick] (11,-0.5) -- (11.5,-0.5);
\node [place,label=below:{\tiny$Q_0^{(0)}$}] at (-0.5,-0.5) {};
\node [place,label=below:{\tiny$P_0^{(0)}$}] at (-1.5,-0.5) {};
\node [place,label=below:{\tiny$P_0^{(1)}$}] at (0.5,-0.5) {};
\node [place,label=below:{\tiny$Q_0^{(1)}$}] at (1.5,-0.5) {};
\node [place,label=below:{\tiny$P_0^{(2)}$}] at (2.5,-0.5) {};
\node [place,label=below:{\tiny$Q_0^{(2)}$}] at (3.5,-0.5) {};
\node [place,label=below:{\tiny$P_0^{(3)}$}] at (4.5,-0.5) {};
\node [place,label=below:{\tiny$\bar{P}_0^{(3)}$}] at (5.5,-0.5) {};
\node [place,label=below:{\tiny$\bar{Q}_0^{(2)}$}] at (6.5,-0.5) {};
\node [place,label=below:{\tiny$\bar{P}_0^{(2)}$}] at (7.5,-0.5) {};
\node [place,label=below:{\tiny$\bar{Q}_0^{(1)}$}] at (8.5,-0.5) {};
\node [place,label=below:{\tiny$\bar{P}_0^{(1)}$}] at (9.5,-0.5) {};
\node [place,label=below:{\tiny$\bar{Q}_0^{(0)}$}] at (10.5,-0.5) {};
\node [place,label=below:{\tiny$\bar{P}_0^{(0)}$}] at (11.5,-0.5) {};
\draw [thick] [->] (-1.5,-0.5)  -- (-1,0);
\draw [thick] (-1,0) --  (-0.5,0.5);
\draw [thick] [->] (-1.5,-0.5)  -- (-0.5,0);
\draw [thick] (-0.5,0) --  (0.5,0.5);
\draw [thick]plot[smooth, tension=.7] coordinates {(-0.5,-0.5) (0.1,-0.1) (0.5,0.5)};
\draw [thick]plot[smooth, tension=.7] coordinates {(-0.5,-0.5) (-0.1,0.1) (0.5,0.5)};
\draw [thick][->] (0.1,-0.1) -- (0.11,-0.09);
\draw [thick][->] (-0.1,0.1) -- (-0.09,0.11);
\draw [thick] [->] (0.5,-0.5) -- (1,0);
\draw [thick] (1,0) --(1.5,0.5);
\draw [thick] [->]  (1.5,-0.5) -- (2,0);
\draw [thick] (2,0) -- (2.5,0.5);
\draw [thick] [->] (2.5,-0.5) -- (3,0);
\draw [thick] (3,0) -- (3.5,0.5);
\draw [thick] [->] (3.5,-0.5) -- (4,0);
\draw [thick] (4,0) -- (4.5,0.5);
\draw [thick] [->] (5.5,0.5)  -- (6,0);
\draw [thick] (6,0) --  (6.5,-0.5);
\draw [thick] [->] (6.5,0.5)  -- (7,0);
\draw [thick] (7,0) --  (7.5,-0.5);
\draw [thick] [->] (7.5,0.5)  -- (8,0);
\draw [thick] (8,0) --  (8.5,-0.5);
\draw [thick] [->] (8.5,0.5)  -- (9,0);
\draw [thick] (9,0) --  (9.5,-0.5);
\draw [thick]plot[smooth, tension=.7] coordinates {(9.5,0.5) (10.1,0.1) (10.5,-0.5)};
\draw [thick]plot[smooth, tension=.7] coordinates {(9.5,0.5) (9.9,-0.1) (10.5,-0.5)};
\draw [thick] [->] (10.09,0.11) -- (10.1,0.1);
\draw [thick] [->] (9.89,-0.09) -- (9.9,-0.1);
\draw [thick] [->] (9.5,0.5)  -- (10.5,0);
\draw [thick] (10.5,0) --  (11.5,-0.5);
\draw [thick] [->] (10.5,0.5)  -- (11,0);
\draw [thick] (11,0) --  (11.5,-0.5);
\end{tikzpicture}
\caption{Digraph $D^{H_3}$}\label{fig-dh3}
\end{figure}
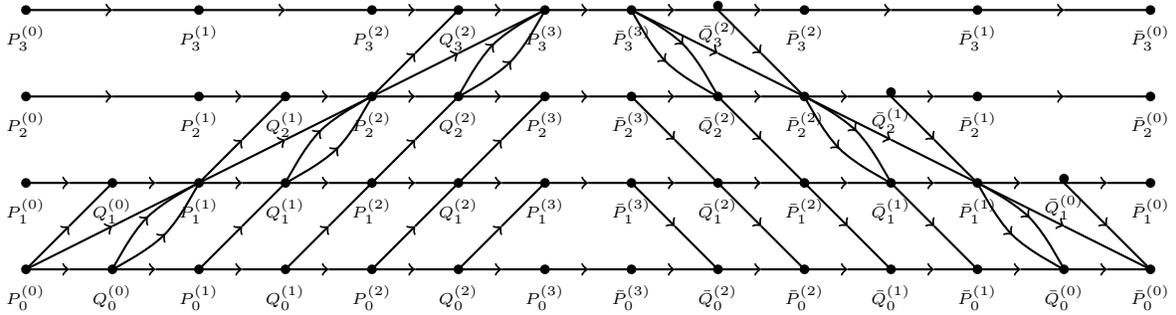

Based on \eqref{eq-Hankel-recurrence} and Lemma \ref{lem-transfer-matrix}, we proceed to build the planar network for $H_n$ from $\mathcal{B}_n$. Firstly, we construct a planar network for $B_n^t$. We take $D^{B_n^t}$ to be the digraph obtained by reflecting $D^{B_n}$ about the vertical line $x = 2n + 1/2$ and reversing the direction of all arcs. We also label the image of $P_j^{(i)}$ (resp. $Q_j^{(i)}$) by $\bar{P}_j^{(i)}$ (resp. $\bar{Q}_j^{(i)}$). We also let $\mathrm{wt}_{D^{B_n^t}}$ be the function which assigns to each arc of $D^{B_n^t}$ the weight of its preimage. Then it is easy to verify that
\[
\mathcal{B}_n^t = (D^{B_n^t}, \mathrm{wt}_{D^{B_n^t}}, (\bar{P}_{n}^{(n)},\ldots,\bar{P}_{0}^{(n)}), (\bar{P}_{n}^{(0)},\ldots,\bar{P}_{0}^{(0)}))
\]
is a planar network for $B_n^t$. Next, we define $D^{T_n}$ to be the digraph whose vertex set is $\{P_{i}^{(n)} \mid 0 \le i \le n\} \cup \{\bar{P}_{j}^{(n)} \mid 0 \le j \le n\}$ and arc set is $\{P_{i}^{(n)} \to \bar{P}_{i}^{(n)} \mid 0 \le i \le n\}$, and let $\mathrm{wt}_{D^{T_n}}(P_{i}^{(n)} \to \bar{P}_{i}^{(n)}) = (T_n)_{i,i}$ for $0 \le i \le n$. Then
\[
\mathcal{T}_n = (D^{T_n},\mathrm{wt}_{D^{T_n}},(P_{n}^{(n)},\ldots,P_{0}^{(n)}),(\bar{P}_{n}^{(n)},\ldots,\bar{P}_{0}^{(n)}))
\]
is a planar network for $T_n$. Finally, we combine $\mathcal{B}_n$, $\mathcal{T}_n$ and $\mathcal{\bar{B}}_n^t$ to get the following planar network for $H_n$:
\begin{align}\label{eq-network-hn}
\mathcal{H}_n = (D^{H_{n}}, \mathrm{wt}_{D^{H_{n}}}, (P_{n}^{(0)},\ldots,P_{0}^{(0)}),(\bar{P}_{n}^{(0)},\ldots,\bar{P}_{0}^{(0)})),
\end{align}
where $D^{H_{n}}$ and $\mathrm{wt}_{D^{H_{n}}}$ are defined in the way as described in Lemma \ref{lem-transfer-matrix}. Figure \ref{fig-dh3} shows the digraph $D^{H_3}$.

We are now in a position to give a combinatorial proof of the $q$-total positivity of $H_n$ for any nonnegative integer $n$.
Given a positive integer $k$ and two sequences $I = (i_1,\ldots,i_k),\,J = (j_1,\ldots,j_k)$ of indices such that $0 \le i_1 < \cdots < i_k\leq n$ and $0 \le j_1 < \cdots < j_k\leq n$, let
\begin{align}\label{eq-pipjbar}
\mathbf{P}_I=(P_{n-i_1}^{(0)},\ldots,P_{n-i_k}^{(0)}),\quad
\mathbf{\bar{P}}_J=(\bar{P}_{n-j_1}^{(0)},\ldots,\bar{P}_{n-j_k}^{(0)}).
\end{align}
Let $H_{I,J}$ denote the submatrix of $H_n$ whose rows are indexed by $I$ and columns indexed by $J$.
By Lemma \ref{lem-lgv} and \eqref{eq-network-hn}, we have
\begin{align}\label{eq-hankel-lgv}
\det \left[H_{I,J}\right] = GF(\mathcal{N}_{D^{H_n}}(\mathbf{P}_I,\mathbf{\bar{P}}_J)).
\end{align}
We further need to find a subset of $\mathcal{N}_{D^{H_n}}(\mathbf{P}_I,\mathbf{\bar{P}}_J)$, say
$\mathcal{S}^{{H_n}}_{I,J}$, which will play the same role as $S_{I,J}$ in Theorem \ref{thm-qtp-lb}. 

Observe that by the recursive construction of $D^{H_n}$, it can be naturally divided into $2n+1$ parts: $D^{H_n}_1$, \ldots, $D^{H_n}_n$, $D^{T_n}$, $\bar{D}^{H_n}_n$, \ldots, $\bar{D}^{H_n}_1$, where for each $1 \le i \le n$ the graph $D^{H_n}_i$ is the subgraph of $D^{H_n}$ induced by the vertices $P_{n}^{(i-1)},\ldots,P_{0}^{(i-1)}$, $Q_{n}^{(i-1)},\ldots,Q_{0}^{(i-1)}$, $P_{n}^{(i)},\ldots,P_{0}^{(i)}$, and $\bar{D}^{H_n}_i$ is the subgraph of $D^{H_n}$ induced by the vertices $\bar{P}_{n}^{(i)},\ldots,\bar{P}_{0}^{(i)}$, $\bar{Q}_{n}^{(i-1)},\ldots,\bar{Q}_{0}^{(i-1)}$, $\bar{P}_{n}^{(i-1)},\ldots,\bar{P}_{0}^{(i-1)}$.
Graphically, $D^{H_n}$ is divided into $2n+1$ parts by $2n$ lines parallel to the $y$-axis. Thus, each member $\mathbf{p}=(p_1,\ldots,p_k)\in \mathcal{N}_{D^{H_n}}(\mathbf{P}_I,\mathbf{\bar{P}}_J)$ is also
divided into $2n+1$ nonintersecting families
$\mathbf{p}_1,\ldots,\mathbf{p}_n, \mathbf{p}_T, \mathbf{\bar{p}}_n,\ldots,\mathbf{\bar{p}}_1$ by these lines,
where $\mathbf{p}_i$ (resp. $\mathbf{\bar{p}}_i$) is the restriction of $\mathbf{p}$ to ${D}^{H_n}_i$ (resp. $\bar{D}^{H_n}_i$) for each
$1 \le i \le n$, and $\mathbf{p}_T$ is the restriction of $\mathbf{p}$ to $D^{T_n}$. For this reason, we may adopt the notation
$$\mathbf{p}=(\mathbf{p}_1,\ldots,\mathbf{p}_n, \mathbf{p}_T, \mathbf{\bar{p}}_n,\ldots,\mathbf{\bar{p}}_1)$$
to represent a nonintersecting family of $\mathcal{N}_{D^{H_n}}(\mathbf{P}_I,\mathbf{\bar{P}}_J)$.

Note that $D^{H_n}_i$ ($1 \le i \le n)$ can be regarded as the digraph obtained by adding $n-i$ parallel arcs (namely, $P_{n}^{(i-1)} \to P_{n}^{(i)},\ldots,P_{i+1}^{(i-1)} \to P_{i+1}^{(i)}$) to $D^{L^{B}_{i-1}}$. By simply mimicking the definitions of ($\mathcal{P}_1$), ($\mathcal{P}_2$) and ($\mathcal{P}_3$) as given immediately before \eqref{eq-pipjp}, we may define the following properties on nonintersecting families $\mathbf{p}_i$ in $D^{H_n}_i$:

\begin{list}{}{\setlength{\leftmargin}{1.1cm}}
\item[($\mathcal{P}_1^{(i)}$)] There exists $1 \le j \le k$ such that the $j$-th component of $\mathbf{p}_i$ is the directed path $P_{i-1}^{(i-1)} \to P_i^{(i)}$;
\item[($\mathcal{P}_2^{(i)}$)] There exists $1 \le j \le k$ such that the $j$-th component of $\mathbf{p}_i$ is the directed path $P_{i-1}^{(i-1)} \to Q_{i-1}^{(i-1)} \overset{l}{\to} P_i^{(i)}$;
\item[($\mathcal{P}^{(i)}_3$)] There exist $1 \le j \le k$ and $l\ge 2$ such that the $(j+l-1)$-th component of $\mathbf{p}_i$ is $P_{i-l}^{(i-1)} \to Q_{i-l}^{(i-1)} \to P_{i-l+1}^{(i)}$ and the $m$-th component is $P_{i-1-(m-j)}^{(i-1)} \to Q_{i-(m-j)}^{(i-1)} \to P_{i-(m-j)}^{(i)}$ for each $j \le m \le j+l-2$.
\end{list}
For $\mathbf{\bar{p}}_i$ in $\bar{D}_i^{H_n}$, if its preimage with respect to the reflection satisfies ($\mathcal{P}^{(i)}_1$), ($\mathcal{P}^{(i)}_2$), or ($\mathcal{P}^{(i)}_3$), we say that $\mathbf{\bar{p}}_i$ satisfies Property ($\mathcal{\bar{P}}^{(i)}_1$), ($\mathcal{\bar{P}}^{(i)}_2$), or ($\mathcal{\bar{P}}^{(i)}_3$), respectively. Then we take
\begin{align}\label{eq-sdhij}
\mathcal{S}^{{H_n}}_{I,J}=\left\{
\mathbf{p}\in \mathcal{N}_{D^{H_n}}(\mathbf{P}_I,\mathbf{\bar{P}}_J) \left |
\begin{array}{l}
\mbox{$\mathbf{p}_i$ satisfies none of ($\mathcal{P}^{(i)}_1$), ($\mathcal{P}^{(i)}_2$), ($\mathcal{P}^{(i)}_3$) and }\\
\mbox{$\mathbf{\bar{p}}_i$ satisfies none of  ($\mathcal{\bar{P}}^{(i)}_1$), ($\mathcal{\bar{P}}^{(i)}_2$), ($\mathcal{\bar{P}}^{(i)}_3$)}\\
\mbox{for each $1\leq i\leq n$}
\end{array}\right.
\right\}.
\end{align}
It is clear that each $\mathbf{p}\in \mathcal{S}^{{H_n}}_{I,J}$
has a $q$-nonnegative weight.

We would like to point out that the involution $\phi$ defined in the proof of Theorem \ref{thm-qtp-lb} can also be mimicked to define a sign-reversing involution $\phi_i$ on nonintersecting families $\mathbf{p}_i$ in $D^{H_n}_i$. Suppose that $\mathbf{p}_i = (p_{i,1},\ldots,p_{i,k})$ and $p_{i,1},\ldots,p_{i,m}$ are those parallel arcs out of $D^{L^{B}_{i-1}}$.
Then $(p_{i,m+1},\ldots,p_{i,k})$ is a nonintersecting family in $D^{L^{B}_{i-1}}$.
If $\phi((p_{i,m+1},\ldots,p_{i,k}))=(p'_{i,m+1},\ldots,p'_{i,k})$, then define
\begin{align}\label{eq-phi-i}
\phi_i(\mathbf{p}_i) = (p_{i,1},\ldots,p_{i,m},p'_{i,m+1},\ldots,p'_{i,k}).
\end{align}
Similarly, we can define a sign-reversing involution $\bar{\phi}_i$ on nonintersecting families $\mathbf{\bar{p}}_i$ in $\bar{D}^{H_n}_i$. Note that if $\mathbf{p}_i$ satisfies Property ($\mathcal{P}^{(i)}_1$), ($\mathcal{P}^{(i)}_2$), or ($\mathcal{P}^{(i)}_3$), then $(p'_{i,m+1},\ldots,p'_{i,k})$ satisfies Property ($\mathcal{P}_1$), ($\mathcal{P}_2$), or ($\mathcal{P}_3$) (with a change of labeling), respectively, and hence $\phi((p_{i,m+1},\ldots,p_{i,k})) \neq (p_{i,m+1},\ldots,p_{i,k})$ and $\phi_i(\mathbf{p}_i) \neq \mathbf{p}_i$. An analogous result holds for $\mathbf{\bar{p}}_i$ and $\bar{\phi}_i$.

The main result of this section is as follows, which provides a combinatorial proof of the $q$-total positivity of $H$.

\begin{thm}\label{thm-main}
Given a nonnegative integer $n$ and two sequences $I = (i_1,\ldots,i_k),\,J = (j_1,\ldots,j_k)$ of indices such that $0 \le i_1 < \cdots < i_k\leq n$ and $0 \le j_1 < \cdots < j_k\leq n$, let $H_{I,J}$ denote the submatrix of $H_n$ whose rows are indexed by $I$ and columns indexed by $J$, let
$\mathbf{P}_I,\mathbf{\bar{P}}_J$ be as given by \eqref{eq-pipjbar}, and let
$\mathcal{S}^{{H_n}}_{I,J}$ be as given by \eqref{eq-sdhij}. Then we have
\begin{align}\label{eq-hankel-lgv-s}
\det \left[H_{I,J}\right] = GF(\mathcal{S}^{{H_n}}_{I,J}),
\end{align}
where $GF(\mathcal{S}^{{H_n}}_{I,J})$ denotes the sum of weights of all elements in $\mathcal{S}^{{H_n}}_{I,J}$.
In particular, $H$ is $q$-totally positive.
\end{thm}

\pf By \eqref{eq-hankel-lgv}, it suffices to give a sign-reversing involution, say $\phi^{H_n}$, on $\mathcal{N}_{D^{H_n}}(\mathbf{P}_I,\mathbf{\bar{P}}_J)$
with $\mathcal{S}^{{H_n}}_{I,J}$ being the set of all fixed points.
We proceed to define $\phi^{H_n}$ by using the aforementioned involutions $\phi_i$ and $\bar{\phi}_i$, see \eqref{eq-phi-i}.
Given $\mathbf{p}=(\mathbf{p}_1,\ldots,\mathbf{p}_n, \mathbf{p}_T, \mathbf{\bar{p}}_n,\ldots,\mathbf{\bar{p}}_1)\in \mathcal{N}_{D^{H_n}}(\mathbf{P}_I,\mathbf{\bar{P}}_J)$, if $\mathbf{p}\in \mathcal{S}^{{H_n}}_{I,J}$, then let $\phi^{H_n}(\mathbf{p})=\mathbf{p}$.

Next, we consider the case $\mathbf{p} \in \mathcal{N}_{D^{H_n}}(\mathbf{P}_I,\mathbf{\bar{P}}_J) \setminus \mathcal{S}^{{H_n}}_{I,J}$. If there exists some $i$ such that $\mathbf{p}_i$ satisfies Property $(\mathcal{P}^{(i)}_1)$, $(\mathcal{P}^{(i)}_2)$ or $(\mathcal{P}^{(i)}_3)$, or equivalently, $\phi_i(\mathbf{p}_i)\neq \mathbf{p}_i$, then let $l$ be the smallest such index and
\[
\phi^{H_n}((\mathbf{p}_1,\ldots,\mathbf{p}_n,\mathbf{p}_T,\bar{\mathbf{p}}_n,\ldots,\bar{\mathbf{p}}_1)) = (\mathbf{p}_1,\ldots,\phi_l(\mathbf{p}_l),\ldots,\mathbf{p}_n,\mathbf{p}_T,\bar{\mathbf{p}}_n,\ldots,\bar{\mathbf{p}}_1).
\]

Otherwise, if such an index does not exist, then there must exist some $i$ such that $\mathbf{\bar{p}}_i$ satisfies Property $(\mathcal{\bar{P}}^{(i)}_1)$, $(\mathcal{\bar{P}}^{(i)}_2)$ or $(\mathcal{\bar{P}}^{(i)}_3)$, or equivalently, $\bar{\phi}_i(\mathbf{\bar{p}}_i)\neq \mathbf{\bar{p}}_i$. In this subcase, we let $l$ be the largest such index and define
\[
\phi^{H_n}((\mathbf{p}_1,\ldots,\mathbf{p}_n,
\mathbf{p}_T,\bar{\mathbf{p}}_n,\ldots,\bar{\mathbf{p}}_1)) = (\mathbf{p}_1,\ldots,\mathbf{p}_n,\mathbf{p}_T,
\bar{\mathbf{p}}_n,\ldots,\bar{{\phi}}_l(\bar{\mathbf{p}}_l),\ldots,\bar{\mathbf{p}}_1).
\]

By the construction of the involutions $\phi_i$ and $\bar{\phi}_i$ for $1 \le i \le n$, it is easy to verify that $\phi^{H_n}$ is a sign-reversing involution on nonintersecting families in $D^{H_n}$. Hence $\phi^{H_n}$ induces the $q$-total positivity of $H_n$ in the same way that $\phi$ induces the $q$-total positivity of $L^B_n$. Further, the $q$-total positivity of $H_n$ implies the $q$-total positivity of $H$ since each minor of $H$ is a minor of $H_n$ for some $n$. The proof is complete. \qed

By applying the same kind of reasoning of the proof of Theorem \ref{thm-main}, we can give a combinatorial proof of the $q$-total positivity of the triangular array
$B=(b_{n,k}(q))_{n,k \ge 0}$ as defined by \eqref{eq-narab-recurrence}. We leave the details to the reader. It is interesting to note that similar reasoning can be used to establish the following result.

\begin{thm}
Let $C = (c_{n,k}(q))_{n,k \ge 0}$ be the Catalan-Stieltjes matrix generated by one of the following two recurrences:
\begin{itemize}
\item[(1)] We have
\begin{align*}
c_{n,0}(q)&= [(f-e)+eq)]\cdot c_{n-1,0}(q) + fq \cdot c_{n-1,1}(q);\\
c_{n,k}(q)&= c_{n-1,k-1}(q) + (1+q)\cdot c_{n-1,k}(q) + q \cdot c_{n-1,k+1}(q), \quad (k\ge 1, \, n\ge 1)
\end{align*}
for some $f \ge e \ge 0$.
\item[(2)] We have
\begin{align*}
c_{n,0}(q)&= [(f-1)+eq)]\cdot c_{n-1,0}(q) + efq \cdot c_{n-1,1}(q);\\
c_{n,k}(q)&= c_{n-1,k-1}(q) + (1+eq)\cdot c_{n-1,k}(q) + eq \cdot c_{n-1,k+1}(q), \quad (k\ge 1, \, n\ge 1)
\end{align*}
for some $e,f \ge 1$.
\end{itemize}
Then the Hankel matrix $(c_{i+j,0})_{i,j \ge 0}$ is $q$-totally positive.
\end{thm}

\pf Note that, for either case, we can construct a planar network for the leading principal submatrix $(c_{i+j,0})_{0 \le i,j \le n}$ by using the
same underlying graph $D^{H_n}$ of the planar network for $H_n$, but with a different weight function. Since the weight function of $\mathcal{H}_n$ is naturally inherited from that of $\mathcal{L^B}$, it is sufficient to assign a new weight function to $D^{L^B}$.

For the first case, we let
\begin{align*}
&\mathrm{wt}_{D^{L^B}}(P_1 \to Q_0) = eq, \quad \mathrm{wt}_{D^{L^B}}(P_i \to Q_{i-1}) = q \text{ for } i \ge 2, \\
&\mathrm{wt}_{D^{L^B}}(P_1 \to P'_0) = -e, \quad \mathrm{wt}_{D^{L^B}}(Q_1 \overset{l}{\to} P'_0) = e, \\
&\mathrm{wt}_{D^{L^B}}(Q_1 \overset{r}{\to} P'_0) = f-e,
\end{align*}
and $\mathrm{wt}_{D^{L^B}}(a) = 1$ for the other arcs $a$ in $D^{L^B}$.

For the second case, we let
\begin{align*}
\mathrm{wt}_{D^{L^B}}(P_i \to Q_{i-1}) = eq \text{ for } i \ge 1, \quad \mathrm{wt}_{D^{L^B}}(P_1 \to P'_0) = -1,\\
\mathrm{wt}_{D^{L^B}}(Q_1 \overset{l}{\to} P'_0) = 1, \quad \mathrm{wt}_{D^{L^B}}(Q_1 \overset{r}{\to} P'_0) = f-1,
\end{align*}
and $\mathrm{wt}_{D^{L^B}}(a) = 1$ for the other arcs $a$ in $D^{L^B}$.

With these new weights, it is straightforward to verify that the involution $\phi^{H_n}$ constructed in the proof of
Theorem \ref{thm-main} is still a sign-reversing involution on nonintersecting families of $D^{H_n}$.
As a result, we obtain the $q$-total positivity of $(c_{i+j,0})_{i,j \ge 0}$. \qed

\section{A conjecture on immanant positivity}
\label{sect-conj}


Let $M = (m_{i,j})$ be a square matrix of order $n$, $\mathfrak{S}_n$ be the symmetric group of order $n$, $\lambda$ be a partition of $n$, and $\chi^{\lambda}$ be the irreducible character of $\mathfrak{S}_n$ associated with $\lambda$. Recall that the immanant of $M$ with respect to $\lambda$ is defined by
\begin{align*}
\mathrm{Imm}_{\,\lambda}\,M = \sum_{\sigma \in \mathfrak{S}_n} \chi^{\lambda}(\sigma) m_{1,\sigma(1)}\cdots m_{n,\sigma(n)}.
\end{align*}
When $\lambda = (1^n)$, the immanant $\mathrm{Imm}_{\,\lambda}\,M$ specializes to $\det [M]$. In \cite{LLYZ21} we proved the immanant positivity for a large family of Catalan-Stieltjes matrices and their associated Hankel matrices. This motivates us to study the
the immanant positivity for the Hankel matrix $H$ defined as in \eqref{eq-Hankel-def}. We have the following conjecture.

\begin{conj}
Let $k \ge 1$ and $I = (i_1,\ldots,i_k)$, $J = (j_1,\ldots,j_k)$ be two sequences of indices with $0 \le i_1 < \cdots < i_k$ and $0 \le j_1 < \cdots < j_k$. Let $H_{I,J}$ be the submatrix of $H$ whose rows are indexed by $I$ and columns are indexed by $J$. Then
\[
\mathrm{Imm}_{\,\lambda}\,H_{I,J} \ge_q 0
\]
for any partition $\lambda$ of $k$.
\end{conj}

We have verified the immanant positivity of all square submatrices of $H_n$ for $n\leq 6$ by Sage \cite{Sage}.
Note that, our method in \cite{LLYZ21} does not apply to $H$ directly, since there exist some arcs weighted by $-1$ in our planar network for $H$.
\vskip 0.5cm
\noindent \textbf{Acknowledgments.} This work is supported in part by the Fundamental Research Funds for the Central Universities and the National Science Foundation of China (Nos. 11522110, 11971249).

\end{document}